\documentclass[12pt]{amsart}
\usepackage{amscd,amssymb}
\usepackage[plainpages,backref,urlcolor=blue]{hyperref}

\theoremstyle{plain}
\newtheorem{thm}[subsection]{Theorem}

\newtheorem{prop}[subsection]{Proposition}

\theoremstyle{definition}
\newtheorem{rk}[subsection]{Remark}

\numberwithin{equation}{section}
\newcommand{\A}{{\mathcal A}}

\newcommand{\al}{{\alpha}}
\newcommand{\be}{{\beta}}

\newcommand{\C}{\mathbb{C}}
\newcommand{\PP}{\mathbb{P}}

\DeclareMathOperator{\mult}{mult}

\begin{document}

\title [On the duals of smooth projective complex hypersurfaces]
{On the duals of smooth projective complex hypersurfaces}

\author[Alexandru Dimca]{Alexandru Dimca$^{1}$}
\address{Universit\'e C\^ ote d'Azur, CNRS, LJAD, France and Simion Stoilow Institute of Mathematics,
P.O. Box 1-764, RO-014700 Bucharest, Romania}
\email{dimca@unice.fr}

\author[Giovanna Ilardi]{Giovanna Ilardi}
\address{Dipartimento Matematica Ed Applicazioni ``R. Caccioppoli''
Universit\`{a} Degli Studi Di Napoli ``Federico II'' Via Cintia -
Complesso Universitario Di Monte S. Angelo 80126 - Napoli - Italia}
\email{giovanna.ilardi@unina.it}


\thanks{\vskip0\baselineskip
\vskip-\baselineskip
\noindent $^1$This work has been partially supported by the Romanian Ministry of Research and Innovation, CNCS - UEFISCDI, grant PN-III-P4-ID-PCE-2020-0029, within PNCDI III}



\begin{abstract}  We  show first that a generic hypersurface $V$ of degree $d\geq 3$ in the complex projective space $ \PP^n$ of dimension $n \geq 3$ has at least one hyperplane section $V \cap H$ containing exactly $n$ ordinary double points, alias $A_1$ singularities, in general position, and no other singularities. Equivalently, the dual hypersurface $V^{\vee}$ has at least one normal crossing singularity of multiplicity $n$. Using this result, we show  that the dual of any smooth hypersurface with $n,d \geq 3$ has at least a very singular point $q$, in particular a point $q$ of multiplicity $\geq n$.

\end{abstract}
 
\maketitle


\section{Introduction} 

Let $S=\C[x_0,\ldots,x_n]$ be the graded polynomial ring in $n+1$ variables with complex coefficients, with $n \geq 2$. Let $f \in S_d$ be a homogeneous polynomial such that the hypersurface $V=V(f):f=0$ in the projective space $\PP^n$ is smooth. Our  main result is the following.
\begin{thm} \label{thm3}
For any $d \geq 3$ and $n \geq 3$, a generic hypersurface $V \subset \PP^n$ of degree $d$ has at least a hyperplane section $V \cap H$, which has exactly $n$ nodal singularities in general position.
\end{thm}
Here a node is a non-degenerate quadratic singularity, in Arnold notation an $A_1$ singularity. A nodal hypersurface is a hypersurface having only such nodes as singularities. As usual, a property is generic if it holds for a Zariski open and dense subset of the parameter space, which is in this case $\PP(S_d)$.
For any smooth hypersurface $V$ and any hyperplane $H$, the singularities of the hyperplane section $V \cap H$ are exactly the points where $H$ is tangent to $V$. The fact that a generic hypersurface has  the tangency points in general position was established in \cite{Bruce}. This means that for a generic, smooth hypersurface $V \subset \PP^n$ and any hyperplane $H \subset \PP^n$, the singularities of the section $V \cap H$ are points in general position in $H$, i.e. the corresponding vectors in the vector space associated to $H$ are linearly independent. In particular, there are at most $n$ singularities in any such section $V \cap H$ when $V$ is generic.

For a smooth hypersurface $V(f)$, consider the corresponding dual mapping
$$\phi_f: V(f) \to (\PP^n)^{\vee}, \  x \mapsto (f_0(x): f_1(x): \ldots : f_n(x)),$$
where we set
$$f_j= \frac{\partial f}{\partial x_j}
\text{ for } j=0, \ldots ,n.$$
Then the dual hypersurface 
$$V(f)^{\vee} = \phi_f (V(f))$$
has a normal crossing singularity of multiplicity $n$ at the point of the dual projective space $(\PP^n)^{\vee}$ corresponding to the hyperplane $H$ if and only if $V(f)\cap H$  has exactly $n$ singularities of type $A_1$ in general position. To prove this, note that $V(f)\cap H$ has a node $A_1$ at a point $p$ if and only if the dual mapping $\phi_f$ is an immersion at $p$,
 see for instance the equivalences $(11.33)$ in \cite{RCS}. Moreover, the tangent space to the corresponding branch $ \phi_f (V(f),p)$ of $V(f)^{\vee}$ at $H=\phi_f (p)$ is given by $$p \in \PP^n=((\PP^n)^{\vee})^{\vee}.$$
 It follows that Theorem \ref{thm3} can be reformulated as follows.

\begin{thm}\label{thm4}
For any dimension $n \geq 3$ and degree $d \geq 3$, the dual hypersurface $V^{\vee}$ of a  generic hypersurface $V \subset \PP^n$ of degree $d$ has at least one normal crossing singularity of multiplicity $n$.
\end{thm}

\begin{rk}\label{rk1}

 (i) For any smooth hypersurface $V$ and any hyperplane $H$,
 the hyperplane section $V\cap H$ has only isolated singularities. Conversely, any hypersurface $W \subset H$ with only isolated singularities may occur as a section 
$W=V \cap H$ for a certain smooth hypersurface $V$, see \cite[Proposition (11.6)]{RCS}.

(ii) The fact that a generic curve $C$ in $\PP^2$ of degree $d\geq 4$ has only simple tangents, bitangents and simple flexes is known classically, and correspond to the claim that the dual curve $C^{\vee}$
has only nodes $A_1$ and cusps $A_2$ as singularities. The number of bitangents of $C$, which is also the number of nodes $A_1$ of the dual curve, as a function of $d$ is also known, see \cite{GH}, page 277 for the classical approach using Pl\"ucker formulas, or \cite{AC} for a modern view-point. It follows from \cite[Proposition 2.1]{Ku}
that {\it any smooth quartic curve has at least 16 bitangents} (which are called there simple bitangents) which correspond to the nodes of the dual curve.

(iii) The fact that a surface $S \subset \PP^3$ of degree $d \geq 3$ admits tritangent planes $H$  is well known, and there are formulas for the number of these planes in terms of the degree $d$, see for instance 
\cite[Section (8.3)]{V}. The fact that, for $S$ generic of degree $d\geq 5$, the singularities of $S \cap H$ are exactly 3 nodes follows from \cite[Proposition 3]{Xu}. For a generic hypersurface $V \subset \PP^n$, with $n \geq 4$ and $d= \deg V \geq n+2$, it follows from
\cite[Proposition 4]{Xu} that all the singularities of a hyperplane section $V \cap H$ are double points (not necessarily $A_1$ singularities) in number {\it at most} $n$.
\end{rk}

Using the above results for generic hypersurfaces, one can prove the following result, {\it holding for any smooth hypersurface}.

\begin{thm}\label{thm5}
For any dimension $n \geq 3$ and degree $d \geq 3$, the dual hypersurface $V^{\vee}$ of a  smooth hypersurface $V \subset \PP^n$ of degree $d$ has either a singularity of multiplicity $n$ with the corresponding tangent cone a union of hyperplanes, or a singularity of multiplicity $>n$. Moreover, 
a  smooth hypersurface $V \subset \PP^n$ of degree $d$, where $n,d \geq 3$, has at least one hyperplane section $V \cap H$ whose total Tjurina number $\tau(V \cap H)$ is at least $n$.
\end{thm}
We recall that for an isolated hypersurface singularity $(X,0):g=0$ defined by a germ $g \in R=\C[[y_1,\ldots,y_n]]$, we define its Milnor number $\mu(X,0)$ and its Tjurina number $\tau(X,0)$ by the formulas
\begin{equation} \label{eqMT} 
\mu(X,0)=\dim R/J_g  \text{ and } \tau(X,0)=\dim R/(J_g+(g)),
\end{equation} 
where $J_g$ is the Jacobian ideal of $g$ in $R$. For a projective hypersurface $W$ having only isolated singularities, we define
 its total Milnor number $\mu(W)$ and its total Tjurina number $\tau(W)$ by the formulas
\begin{equation} \label{eqMT2} 
\mu(W)=\sum_p\mu(W,p)  \text{ and } \mu(W)=\sum_p\mu(W,p),
\end{equation} 
where both sums are over all the singular points $p \in W$. For any point $H \in V^{\vee}$ it is known that
\begin{equation} \label{eqMT3} 
\mult_H(V^{\vee})=\mu(V \cap H),
\end{equation} 
where $\mult_p(Y)$ denotes the multiplicity of a variety $Y$ at a point $p \in Y$, see \cite{DD}. Hence the first claim in Theorem \ref{thm5} implies that 
$$\mu (V \cap H) \geq n$$
 if $\mult_H(V^{\vee}) \geq n$. However,
our second claim in Theorem \ref{thm5} is a stronger version of this inequality, since $\mu(X,0) \geq \tau(X,0)$, with equality exactly when the singularity $(X,0)$ is weighted homogeneous, see \cite{KS}.

We would  like to thank the referee for the very careful reading of our manuscript and for his suggestions to improve the presentation.

\section{Proof of Theorem \ref{thm3}} 
The starting point is Remark \ref{rk1} (i) above.
We consider first the projective space $\PP^{n-1}$ and the subset
$Z_n \subset \PP^{n-1}=\PP(\C^n)$ given by the classes $p_i$ of the canonical basis $e_i$, $i=1, \ldots ,n$ of the vector space $\C^n$.

\begin{prop}\label{prop21}
For any degree $d \geq 3$, there is a hypersurface $Y \subset \PP^{n-1}$, with $n \geq 3$,  of degree $d$ having as singularities $n$ nodes $A_1$, located at the points in $Z_n$.
\end{prop}
\proof Let $y=(y_1, \ldots,y_n)$ be the coordinates on $\PP^{n-1}$.
We consider first the case $d=3$ and take $Y$ to be the hypersurface $g(y)=0$, where
$$g(y)= \sum_{1 \leq i <j <k\leq n}y_iy_jy_k.$$
It is easy to see that $Y$ has an $A_1$-singularity at each point $p_i$, 
for $i=1, \ldots ,n$. Now we show that there are no other singularities.
Note that for the partial derivative $g_i$ of $g$ with respect to $y_i$ we have 
\begin{equation} \label{e21} 
g_i(y)= \sum_{1 \leq j <k\leq n, j \ne i, k \ne i}y_jy_k.
\end{equation} 
Assume that $g_i(y)=0$ for $i=1, \ldots ,n$ and take the sum of all these equations. We get in this way
\begin{equation} \label{e22} 
\sum_{1 \leq j <k\leq n}y_jy_k=0.
\end{equation} 
Subtracting the equation \eqref{e21} from \eqref{e22} we get
\begin{equation} \label{e23} 
y_i\sum_{1 \leq j \leq n, j \ne i}y_j=0.
\end{equation} 
If we assume that $y_{i_1} \ne 0$ and $y_{i_2} \ne 0$ for some indices
$1 \leq i_1<i_2 \leq n$, the equation \eqref{e23} implies that $y_{i_1}=y_{i_2}$. Hence, for any singular point $y^0$ of $Y$, there is an integer
$a$ with $1 \leq a \leq n$ such that $a$ coordinates of $y^0$ are equal to 1, and the remaining $n-a$ coordinates are 0. The equation  \eqref{e23} implies that only the case $a=1$ is possible, and hence $y^0$ is one of the points $p_i$. This completes the proof in the case $d=3$.

Next we look at the case $d=4$  and take $Y$ to be the hypersurface $g(y)=0$, where
$$g(y)= \sum_{1 \leq i <j \leq n}y_i^2y_j^2.$$
It is easy to see that $Y$ has an $A_1$-singularity at each point $p_i$, 
for $i=1, \ldots ,n$. Now we show that there are no other singularities.
In this case we have
\begin{equation} \label{e24} 
g_i(y)= 2y_i\sum_{1 \leq j \leq n, j \ne i}y_j^2.
\end{equation} 
If we assume $g_i(y)=0$ for all $i$, and that $y_{i_1} \ne 0$ and $y_{i_2} \ne 0$ for some indices
$1 \leq i_1<i_2 \leq n$, the equation \eqref{e24} implies that $y_{i_1}^2=y_{i_2}^2$. Hence, for any singular point $y^0$ of $Y$, there is an integer
$a$ with $1 \leq a \leq n$ such that $a$ coordinates of $y^0$ are equal to $\pm 1$, and the remaining $n-a$ coordinates are 0. The equation  \eqref{e24} implies that only the case $a=1$ is possible, and hence $y^0$ is one of the points $p_i$. This completes the proof in the case $d=4$.

Finally, to treat the case $d>4$, let 
$$h_i(y)=y_i^{d-2} \sum_{1 \leq j \leq n, j \ne i}y_j^2.$$
Note that $h_i$ has a singularity of type $A_1$ at $p_i$ and vanishes of order $d-2>2$ at the other points $p_j$, for $j \ne i$.
Consider the linear system spanned by $h_1, \ldots, h_n$. It is easy to see, repeating the  argument already used twice above, that the base locus $h_1= \ldots =h_n=0$ of this linear system is exactly the set $Z_n$. It follows, by Bertini's Theorem, that a generic member $Y$ of this linear system is smooth except possibly at the points of $Z_n$.
The choice of the $h_i$ implies that $Y$ has an $A_1$ singularity at each point in $Z_n$.

\endproof

Now we give a proof of Theorem \ref{thm3} stated in the Introduction.
 Using Remark \ref{rk1} (i) and Proposition \ref{prop21}, it follows that, in any dimension $n \geq 3$ and degree $d \geq 3$, there are smooth hypersurfaces $V(f) \subset \PP^n$ of degree $d$ which have at least one nodal hyperplane section $V(f) \cap H$, with exactly $n$ singularities in general position. 
 
 Let $B=\PP(S_d)_0$ be the set of points in $\PP(S_d)$ corresponding to  polynomials $f \in S_d$ such that the hypersurface $V(f) :f=0$ is  smooth.
Let $\A(n,d) \subset B$ be the subset of such hypersurfaces $V(f)$, which have at least one hyperplane section $V(f) \cap H$, with exactly $n$ singularities $A_1$  in general position. We know already that $\A(n,d) \ne \emptyset$.  It is easy to see that $\A(n,d)$ is a constructible (or semialgebraic) subset in $B$. Indeed, consider the subset
$$\Gamma \subset B \times (\PP^n)^n$$
consisting of pairs $(f,q)$, where $f \in B$, $q=(q^1, \ldots, q^n) \in (\PP^n)^n$ such that the points 
$$q^j=(q^j_0:q^j_1:\ldots:q^j_n)\in \PP^n$$ 
are linearly independent, that is they span a hyperplane $H(q)$ in $\PP^n$, and the following conditions hold
\begin{equation} \label{E1} 
\sum_{i=0}^nq^j_if_i(q^k)=0 \text{ for any } j,k=1, \ldots, n
\end{equation} 
and 
\begin{equation} \label{E2} 
Hess(f)(q^j)\ne 0 \text{ for any } j=1, \ldots, n,
\end{equation} 
where $Hess(f)$ is the Hessian polynomial of $f$.
In fact, the equation \eqref{E1} for $k=j$ tells us that $q^j \in V(f)$ for any $j=1, \ldots,n$. Moreover, it says that the point $q^j$ is on the tangent space $T_{q^k}V(f)$. This implies that 
\begin{equation} \label{E3} 
T_{q^k}V(f)=H(q) \text{ for any } k=1, \ldots, n.
\end{equation} 
The equation \eqref{E2} tells us that the singularity of $V(f) \cap H(q)$ at the point $q^j$ is a node, see for instance \cite[Equivalence (11.33)]{RCS}. It is clear that $\Gamma$ is a constructible set in $B \times (\PP^n)^n$, since it is defined by finitely many algebraic equalities and inequalities. Let $p_1:B \times (\PP^n)^n \to B$ be the first projection and note that $\A(n,d)=p_1(\Gamma)$. Using Chevalley Theorem, see for instance \cite{Lo}, p. 395, we conclude that the set $\A(n,d)$ is constructible in $B$.

We  show now that $\A(n,d)$ is a non empty Zariski open subset in $B$, and hence it is dense in $B$ and in $\PP(S)_d$, see \cite[Theorem 2.33]{Mu}. Let $Z=B \setminus \A(n,d)$. Then $Z$ is also a constructible set, and  \cite[Proposition 2]{Lo} page 394, implies that the closure of
$Z$ in $B$ in the Zariski topology coincides with its closure in the strong complex topology. Hence, to show that $\A(n,d)$ is a  Zariski open subset in $B$, it is enough to show that $\A(n,d)$ is 
 an open subset of $B$ in the strong complex topology.

We fix now one element $f \in \A(n,d)$ and show that $\A(n,d)$ contains a neighborhood of $f$ in $B$. 
The set $B$ is open, hence there are arbitrarily small open neighborhoods $U$ of $f$  with $f \in U \subset B$. For any polynomial $f' \in U$, we consider the gradient map
$$\Phi_{f'}:\PP^n \to (\PP^n)^{\vee} \text{ given by } x  \mapsto (f_0'(x): f_1'(x): \ldots : f_n'(x))$$
and the corresponding dual mapping
$$\phi_{f'}=\Phi_{f'}|V(f'): V(f') \to (\PP^n)^{\vee}.$$
Since $f \in \A(n,d)$, there is a hyperplane $H$ such that $V(f) \cap H$ has $n$ singularities $A_1$ in general position, say at the points $p_j\in \PP^n$, for $j=1,\ldots,n$. It follows that the Hessian polynomial $Hess(f)$ of $f$ satisfies $Hess(f)(p_j) \ne 0$ for $j=1,\ldots,n$ and hence each analytic germ 
$$\Phi_f:(\PP^n,p_j) \to ((\PP^n)^{\vee},H)$$
is invertible, see for instance \cite[ Equation (11.10)]{RCS}.
It follows that there is a neighborhood $N$ of $H$ in $(\PP^n)^{\vee}$ and neighborhoods $N_j$ of $p_j$ in $\PP^n$ for $j=1,\ldots,n$ such that the restrictions
$$\Phi_f^{(j)}=\Phi_f | N_j :N_j \to N$$
 are analytic isomorphisms,  with corresponding inverse mappings
 $$\Psi_f^{(j)}=(\Phi_f^{(j)})^{-1}:N \to N_j.$$
Any polynomial $f' \in U$ can be regarded as a deformation of $f$, the parameters being the coefficients of $f'$. Since $f$ depends analytically
on these parameters, the inverse mapping $\Psi_f^{(j)}$ also depends analytically on these parameters. It follows that, by choosing small enough neighborhoods $U$, $N$ and $N_j$ for $j=1,\ldots,n$, we have
inverse mappings as above
$$\Psi_{f'}^{(j)}:N \to N_j$$
for all $f' \in U$. Choose $g \in S_d$ a polynomial such that $g(p_j) \ne 0$ for $j=1,\ldots,n$. Then
$$h_{f'}=\frac{f'}{g}$$
is an analytic function defined on all $N_j$'s, if they are chosen small enough. Define now
$$\al_{f'}:N \to \C^n \text{ given by } y  \mapsto (h_{f'}(\Psi_{f'}^{(1)}(y)), \ldots, h_{f'}(\Psi_{f'}^{(n)}(y))).$$
Notice that we have an obvious equality of (possibly non reduced) analytic spaces $\al_f^{-1}(0)=\{H\}$, where $H$ is regarded as a point with its reduced structure.
Indeed, $h_{f}=0$ in $N_j$ defines the intersection $V(f) \cap N_j$, and
$\Phi_{f}^{(j)}(V(f) \cap N_j)$ is the trace on $N$ of the irreducible smooth branch of the dual variety $V(f)^{\vee} $ at the point $H\in  (\PP^n)^{\vee}$, whose tangent space at $H=\Phi_{f}^{(j)}(p_j) $ corresponds to the point $p_j\in \PP^n$.
The intersection of these $n$ smooth branches, meeting transversally at $H$, is exactly the simple point $H$. Let $D$ be a small closed ball in $N$, centered at $H$ and consider the restricted mapping
$$\be_{f'}=\al_{f'}| \partial D: \partial D=S^{2n-1} \to \C^n \setminus \{0\}.$$
Here $\partial D$ is the boundary of the closed ball $D$. If $D$ is small, it is clear by the above discussion that $\al_f$ has no zeros on the compact set $\partial D$. By continuity, the same is true for $\al_{f'}$, and hence $\be_{f'}$ is correctly defined for $f'$ close to $f$.
Notice that the mapping $\be_f$ has degree one, see for instance \cite[Section 5.4]{AGV1} or the topological interpretation of intersection multiplicity of $n$ divisors in \cite{GH}, p. 670.
By continuity, it follows that $\deg \be_{f'}=1$ for any $f' \in U$. Therefore, for any $f' \in U$ there is a unique point $H' \in  (\PP^n)^{\vee}$ such that $\al_{f'}^{-1}(0)=\{H'\}$. As explained above, this is equivalent to the fact that the dual hypersurface  $V(f')^{\vee} $ has
 a normal crossing singularity of multiplicity $n$ at the point $H'$. Hence $H'$ corresponds to a hyperplane section of $V(f')$ with $n$ nodes, that is 
$f' \in \A(n,d)$. Therefore $U \subset  \A(n,d)$ and this completes our proof.
\endproof

\section{Proof of Theorem \ref{thm5}} 
Fix $V(f):f=0$ a smooth hypersurface of degree $d$ in $\PP^n$, with $d,n \geq 3$. In view of Theorem \ref{thm3}, there is a sequence of pairs
$(V(f_m), H_m)$ such that $f_m$ converges to $f$ in the projective space $\PP(S_d)$ and $V(f_m) \cap H_m$ has $n$ nodes in general position. By passing to a subsequence, we can assume that $q_m=H_m$ converges to a hyperplane $q=H$ in $(\PP^n)^{\vee} $.
By passing to the dual hypersurfaces, we get a sequence of hypersurfaces $V(f_m)^{\vee}$ converging to $V(f)^{\vee}$, and a sequence of points $q_m \in V(f_m)^{\vee}$ converging to the point
$q \in V(f)^{\vee}$ such that $(V(f_m)^{\vee},q_m)$ is a normal crossing singularity of multiplicity $n$. If we consider the $n$-jet at $q_m$ of a reduced defining equation $F_m$ for $V(f_m)^{\vee}$ in $(\PP^n)^{\vee} $, we see that
$$h_m=j^n_{q_m}F_m$$
is a degree $n$ homogeneous polynomial which splits as a product of $n$ linearly independent linear forms. Let $F$ be a reduced defining equation  for $V(f)^{\vee}$ in $(\PP^n)^{\vee} $. Then $F_m$ converges to $F$ in $\PP(S_D)$, where $D=d(d-1)^{n-1}$, see for instance \cite[Theorem 1.2.5]{Do}.
Hence for the $n$-jet
$$h=j^n_{q}F$$
there are the following two possibilities.
Either $h \ne 0$, and then $h=\lim h_m$ in the corresponding projective space, and so $h$ is a product of $n$ (maybe non-distinct) linear forms, or $h=0$ and then $\mult_qV(f)^{\vee} \geq n+1$. This proves the first claim in Theorem \ref{thm5}.
To prove the second claim, note that for $m$ large, we can identify the hyperplane $H_m$ with $H$ using a linear projection, and in this way
$V(f_m) \cap H_m$ give rise to a sequence of hypersurfaces $W_m$ in $H$ converging to the intersection $W=V(f) \cap H$.
Let $G_m$ (resp. $G$) be the reduced defining equation of the hypersurface $W_m$ (resp. $W$) in $H= \PP^{n-1}$.
If we choose a system of coordinates $y=(y_1: \ldots:y_n)$ and set
$S'=\C[y_1,\ldots,y_n]$, let $M(G_m)=S'/J(G_m)$ and $M(G)=S'/J(G)$ denote the corresponding Milnor (or Jacobian) algebras of $G_m$ and $G$. Here $J(G_m)$ (resp. $J(G)$) denotes the Jacobian ideal of $G_m$ (resp. $G$) spanned by all the first order partial derivatives of $G_m$ (resp. $G$) with respect to the $y_j$'s. For $k>0$ an integer, we set
$$M(G_m)^k=\frac{S'}{J(G_m)+M^{k+1}} \text{ and }M(G)k=\frac{S'}{J(G)+M^{k+1}},$$
where $M$ is the maximal ideal $(y_1,\ldots,y_n) \subset S'$.
Let $T=n(d-2)$ and recall that the homogeneous components of the Milnor algebras $M(G_m)$ and $M(G)$ satisfy
$$\dim M(G_m)_j= \tau(W_m) \text{ and }\dim M(G)_j= \tau(W),$$
for any $j >T$, see \cite[Corollary 9]{CD}. It follows that
$$\dim M(G_m)^k=\dim M(G_m)^T+(k-T)\tau(W_m)$$
and
$$\dim M(G)^k=\dim M(G)^T+(k-T)\tau(W),$$
for any $k>T$. Using the semicontinuity of the dimension of a quotient space, we get
$$\dim M(G)^k \geq \dim M(G_m)^k,$$
for all $k>T$. This clearly implies
$$\tau(W) \geq \tau(W_m) =n$$
and this proves our second claim.

\end{document}